# Harmonic mappings into non-negatively curved manifolds

**Sergey Stepanov • Irina Tsyganok**

**Abstract** So far, all known results on harmonic maps between Riemannian manifolds are based in an essential way on the assumption that the target manifold has non-positive sectional curvature. In our paper we develop a theory of harmonic mappings of Riemannian manifolds into non-negatively curved Riemannian manifolds and give the geometric applications of our results to the theory of holomorphic maps of almost Kählerian manifolds and to addressing the "prescribed Ricci curvature problem".

## 1. Introduction and results

In this section we focus our attention on the applications of the *Bochner technique* to harmonic maps. We also announce our results on harmonic maps that we have obtained using the Bochner technique which are published in the present paper.

The Bochner technique is the most important analytic method of differential geometry "in the large" which derived by Bochner for proving so-called *vanishing theorems* under appropriate curvature conditions on compact Riemannian manifolds (see [25]; [32] and [33]). In the works of Lichnerowicz, Nomizu, Kodaira, Yano and others, the analytic Bochner method was substantially developed and successfully applied to complex, complete Riemannian, and Lorentzian manifolds. In particular, the Bochner technique has applications in the theory of harmonic maps "in the large" which presented in the monographs [13] and [18] from the present point of view. At the same time, we note that the reader is referred to [1, pp. 65-105]; [5]; [6] and [7] for a detailed account of harmonic maps.

---

S. Stepanov[1,2]

[1] Department of Mathematics, All Russian Institute for Scientific and Technical Information of the Russian Academy of Sciences,
   20, Usievicha street, 125190 Moscow, Russian Federation

[2] Department of Mathematics, Finance University, Leningradsky Prospect, 49-55, 125468 Moscow, Russian Federation
   e-mail: s.e.stepanov@mail.ru

Tsyganok Irina
Department of Mathematics, Finance University, Leningradsky Prospect, 49-55, 125468 Moscow, Russian Federation
e-mail: *i.i.tsyganok@mail.ru*

The first application of the Bochner technique to study harmonic maps was in the well known paper of Eells and Sampson [7] where they was obtained the *Weitzenböck formula* for harmonic maps and proved the celebrated vanishing theorem on harmonic maps which state the following: if $f:(M,g)\to(\overline{M},\overline{g})$ is any harmonic mapping between a compact Riemannian manifold (*M, g*) with the Ricci tensor $Ric \geq 0$ and a Riemannian manifold $(\overline{M},\overline{g})$ with the sectional curvature $\overline{sec} \leq 0$ then *f* is *totally geodesic* (see for the definition [31]) and has constant the *energy density* $e(f)$. Furthermore, if there is at least one point of *M* at which its Ricci curvature $Ric > 0$, then every harmonic map $f:(M,g)\to(\overline{M},\overline{g})$ is constant (see [7]). The most recent vanishing theorem was proved in [16]. It states the following: let (*M, g*) and $(\overline{M},\overline{g})$ be compact Riemannian manifolds with *sec* > 0 and $\overline{sec} \leq 0$, then any map (in particular, harmonic map) from (*M, g*) to $(\overline{M},\overline{g})$ must be homotopic to a constant map.

The above scheme can be extended to a harmonic mapping of a complete manifold to a compact manifold with the non-positive sectional curvature. For example, Yau and Schoen showed the following vanishing theorem (see [22]): A harmonic map of finite energy $E(f)$ from a complete non-compact manifold (*M, g*) with the Ricci tensor $Ric \geq 0$ to a compact manifold $(\overline{M},\overline{g})$ with the sectional curvature $\overline{sec} \leq 0$ is homotopic to a constant map on each compact set.

To summarize the above observations (see also [13] and [18]) we note that until today all known results on harmonic maps between Riemannian manifolds are based in an essential way on the assumption that the target manifold $(\overline{M},\overline{g})$ has *non-positive* sectional curvature.

To contrast the above mentioned in the third section of the present paper we prove that any harmonic mapping $f:(M,g)\to(\overline{M},\overline{g})$ between Riemannian manifolds (*M, g*) and $(\overline{M},\overline{g})$ is totally geodesic if the section curvature of $(\overline{M},\overline{g})$ is *non-negative* and (*M, g*) is a compact manifold with the Ricci tensor $Ric \geq f^*\overline{Ric}$ for the pullback $f^*\overline{Ric}$ of the Ricci tensor $\overline{Ric}$ by *f*. Beside it, in the fifth section we will show that a harmonic mapping of finite energy $E(f)$ from of a complete manifold $(M, g)$ with the Ricci tensor $Ric \geq f^*\overline{Ric}$ to a manifold $(\overline{M}, \overline{g})$ with the sectional curvature $\overline{sec} \geq 0$ is totally geodesic.

The map $f:(M,g) \to (\overline{M},\overline{g})$ is called a *contraction* or *weakly length decreasing* if $f^*\overline{g} \leq g$ everywhere on $M$ (see [21]). In particular, the map $f$ is called *strictly length decreasing* if $f^*\overline{g} < g$. In [16] was proved the following vanishing theorem: Any strictly length decreasing map between compact Riemannian manifolds $(M, g)$ and $(\overline{M}, \overline{g})$ whose sectional curvatures are bounded by $sec \geq \sigma \geq \overline{sec}$ for some positive number $\sigma > 0$, is homotopic to a constant map. Other vanishing theorems can be found in [21] and [30]. In the fourth section of our paper we prove that any harmonic contraction mapping between a compact Riemannian manifold $(M,g)$ with the Ricci tensor $Ric > (\overline{n}-1)g$ and a unite sphere $\mathbb{S}^{\overline{n}}(1)$ is a constant mapping. In particular, if $f: \mathbb{S}^n(1) \to \mathbb{S}^{\overline{n}}(1)$ is a harmonic strictly length decreasing mapping for the case $n > \overline{n}$, then it is constant.

In the case of *Kähler manifolds* Sealey in [23] studied *holomorphic maps* (which are harmonic) between Kähler manifolds $(M, g, J)$ and $(\overline{M}, \overline{g}, \overline{J})$. Another generalization of the Bochner technique was presented by Siu in [24] where he has obtained a Weitzenböck formula involving only the curvature of the image manifold for harmonic maps between Kähler manifolds $(M, g, J)$ and $(\overline{M}, \overline{g}, \overline{J})$. This type of argument enabled him to study properties of harmonic mappings and to obtain rigidity results between compact Kähler manifolds with curvature conditions on the image manifold only. A similar analysis was carried out by Sampson in [20] for harmonic maps from a compact Kähler manifold $(M, g, J)$ into a Riemannian manifold $(\overline{M},\overline{g})$. In addition, we state that the reader is referred to [13] and [18] for a detailed account of harmonic maps between Kähler manifolds and their generalizations.

In the sixth section of our paper we prove that any holomorphic mapping $f:(M,g,J) \to (\overline{M},\overline{g},\overline{J})$ between a complete *almost semi-Kählerian manifold* (M, g,

$J$) with the Ricci tensor $Ric \geq f^*\overline{Ric}$ and a *nearly-Kählerian manifold* $(\overline{M}, \overline{g}, \overline{J})$ with the sectional curvature $\overline{sec} \geq 0$ is totally geodesic if its energy $E(f)$ is finite.

The point of the famous papers [4] and [12] and the monograph [2; pp. 140-153] is that in certain circumstances the metric (or at last the connection) is uniquely determined by the Ricci tensor. In particular, was proved some vanishing theorems. In the last section of the our paper we consider a compact Riemannian manifold $(M, \overline{g})$ with the sectional curvature $\overline{sec} \geq 0$ and the Ricci tensor $\overline{Ric} \leq \overline{g}$. Under these conditions, we prove that if $g$ is another Riemannian metric on $M$ with $Ric = \overline{g}$, then $g$ and $\overline{g}$ have the same Levi-Civita connection. Furthermore, if the full holonomy group $\text{Hol}(\overline{g})$ is irreducible then $\overline{Ric} = Ric$. This proposition generalizes one of the main results of [13].

In the present paper we continue the study which we began in [26] where we proved some vanishing theorems for harmonic, umbilical and projective mappings. At the same time, this paper is based on our report [29] at the 12th International Conference on Geometry and Applications (September 1-5, 2015, Varna, Bulgaria).

**2. Definitions and notations**

In this section we give a brief survey of differentiable mappings of Riemannian manifolds (see, for example, [5]; [14; pp. 8-10]; [18]).

Let $(M, g)$ be a Riemannian manifold of dimensional $n$ with the Levi-Civita connection $\nabla$ and $(\overline{M}, \overline{g})$ be a Riemannian manifold of dimensional $\overline{n}$ with the Levi-Civita connection $\overline{\nabla}$. Let there be given $f:(M,g) \to (\overline{M}, \overline{g})$ a differentiable mapping between Riemannian manifolds. Suppose $f^*T\overline{M}$ is the vector bundle over $M$ with fiber $T_{f(x)}\overline{M}$ over $x \in M$. Then the differential $f_*:TM \to T\overline{M}$ of the mapping $f$ is a smooth section of the bundle $T^*M \otimes f^*T\overline{M}$ over $M$. On the other hand, the transpose of $(f_*)_x$ is a linear mapping of $T^*_{f(x)}\overline{M}$ into $T^*_xM$ at each point $x \in M$. In this case, for any covariant tensor $\overline{\omega}$ on $\overline{M}$, we can define a tensor $f^*\overline{\omega}$ on $M$ by $(f^*\overline{\omega})(X,...,Y) = \overline{\omega}(f_*X,...,f_*Y)$ for an arbitrary $X,...,Y \in T_xM$.

Take coordinate neighborhoods $U \subset M$ with local coordinates $x^1,...,x^n$ and $\overline{U} \subset \overline{M}$ with local coordinates $\overline{x}^1,...,\overline{x}^{\overline{n}}$ such that $f(U) \subseteq \overline{U}$. We denote by $g_{ij}$ the local components of the Riemannian metric $g$ on $U \subset M$, by $\overline{g}_{\alpha\beta}$ those of the Riemannian metric $\overline{g}$ on $\overline{U} \subset \overline{M}$ where the indices $i, j, k, l, \ldots$ run over the range $\{1, \ldots, n\}$ and $\alpha, \beta, \gamma, \delta,...$ run over the range $\{1, \ldots, \overline{n}\}$.

Suppose that $f:(M,g) \to (\overline{M},\overline{g})$ is representation by equations $\overline{x}^\alpha = f^\alpha(x^1,...,x^n)$ with respect to the local coordinates of $U$ and $\overline{U}$. We put $f_i^\alpha = d\overline{x}^\alpha/dx^i$ then the differential $f_*$ of mapping $f$ is represented by the matrix $(f_i^\alpha)$ with respect to the local coordinates of $U$ and $\overline{U}$. Then the *energy density* of $f:(M,g) \to (\overline{M},\overline{g})$ is the non-negative scalar function $e(f): M \to \mathbb{R}$ such that $e(f) = \frac{1}{2} \|f_*\|^2$ where $\|f_*\|^2$ denotes the squared norm of the differential $f_*$, with respect to the induced metric on the vector bundle $T^*M \otimes f^*T\overline{M}$. In local coordinates $e(f)$ is expressed by $e(f) = \frac{1}{2} \|f_*\|^2 = \frac{1}{2} g^{ij} \overline{g}_{\alpha\beta} f_i^\alpha f_j^\beta$ for the local components $g^{ij}$ of the tensor $g^{-1}$.

Let $\Omega \subset M$ be a compact domain then we define the *energy of the map* $f|_\Omega : (\Omega,g) \to (\overline{M},\overline{g})$ by the equality $E_\Omega(f) = \int_\Omega e(f)\, d\operatorname{Vol}_g$ for the canonical measure $d\operatorname{Vol}_g = \sqrt{\det g}\, dx^1 \wedge ... \wedge dx^n$ which associated to $g$ (see [18; p. 1]). A smooth mapping $f:(M,g) \to (\overline{M},\overline{g})$ is said to be *harmonic* if, for each compact domain $\Omega \subset M$, it is a stationary point of the energy functional $E_\Omega : C^\infty(M, \overline{M}) \to \mathbb{R}$ with respect to variations preserving $f$ on $\partial\Omega$ (see [6; p. 389] and [18; p. 1]). It is well known that $f:(M,g) \to (\overline{M},\overline{g})$ is a harmonic mapping if and only if it satisfies the *Euler-Lagrange equation* $\operatorname{trace}_g Df_* = 0$ where $D$ is the connection in the bundle $T^*M \otimes f^*T\overline{M}$ induced from the Levi-Civita connections $\nabla$ and $\overline{\nabla}$ of $(M,g)$ and $(\overline{M},\overline{g})$, respectively (see [7, p. 116]; [6; p. 389]).

# 3. Harmonic mappings of compact Riemannian manifolds to non-negatively curved Riemannian manifolds

In this section, we will show the main theorem of our work, which we have proved by the Bochner technique.

Let $f:(M,g)\to(\overline{M},\overline{g})$ be a harmonic mapping between Riemannian manifolds and assume that $(M,g)$ to be compact. We recall the well known the Eells-Sampson equation for a harmonic mapping $f$ is the following (see [7; p. 123]; [18; p. 3]):

$$\Delta e(f) = Q(f) + \|Df_*\|^2 \qquad (3.1)$$

where $\Delta$ is the Laplace–Beltrami operator $\Delta = div\,\nabla$ and

$$Q(f) = g^{ik}g^{jl}\left(-(f^*\overline{R})_{ijkl} + R_{ij}(f^*\overline{g})_{kl}\right)$$
$$= -g^{ik}g^{jl}\overline{R}_{\alpha\beta\gamma\delta}f_i^\alpha f_j^\beta f_k^\gamma f_l^\delta + g^{ik}g^{jl}R_{ij}f_k^\alpha f_l^\beta \overline{g}_{\alpha\beta} \qquad (3.2)$$

for the local components $\overline{R}_{\alpha\beta\gamma\delta}$ of the Riemannian curvature tensor $\overline{R}$ of $(\overline{M},\overline{g})$ and the local components $R_{ij}$ of the Ricci tensor $Ric$ of $(M,g)$.

Now, using (3.1) and (3.2) we shall prove the vanishing theorem for harmonic maps which is an analogue of the Eells and Sampson celebrated vanishing theorem (see [7; p. 124]; [18; p. 3]). First, we suppose that $\overline{n} \geq 2$ and rewrite (3.2) in the form

$$Q(f) = -\overline{R}_{\alpha\beta\gamma\delta}\Phi^{\alpha\gamma}\Phi^{\beta\delta} + R_{ij}f_k^\alpha f_j^\beta g^{ik}g^{jl}\overline{g}_{\alpha\beta} \qquad (3.3)$$

where $\Phi^{\alpha\gamma} = f_i^\alpha f_k^\gamma g^{ik} = \Phi^{\gamma\alpha}$. Second, we diagonalize the symmetric tensor $\Phi$ with respect to $\overline{g}$, using an orthonormal basis $\{\overline{e}_1,...,\overline{e}_{\overline{n}}\}$ at each point $f(x)\in\overline{M}$. Then letting sectional curvature of the plane $\pi$ of $T_{f(x)}\overline{M}$ generate by $\overline{e}_\alpha$ and $\overline{e}_\beta$ be $\overline{sec}(\pi) = \overline{sec}(\overline{e}_\alpha,\overline{e}_\beta)$ and expressing $\overline{sec}(\overline{e}_\alpha,\overline{e}_\beta)$ in term of the curvature tensor $\overline{R}$ by [2, p. 436] and [3] we have

$$\overline{R}_{\alpha\beta}\overline{g}_{\gamma\delta}\Phi^{\alpha\gamma}\Phi^{\beta\delta} - \overline{R}_{\alpha\beta\gamma\delta}\Phi^{\alpha\gamma}\Phi^{\beta\delta} = \sum_{\alpha<\beta}\overline{sec}(\overline{e}_\alpha,\overline{e}_\beta)(\overline{\lambda}_\alpha - \overline{\lambda}_\beta)^2$$

where all eigen-values $\overline{\lambda}_\alpha \geq 0$. In this case, the equality (3.3) can be rewritten in the form

$$Q(f) = \sum_{\alpha<\beta}\overline{sec}(\overline{e}_\alpha,\overline{e}_\beta)(\overline{\lambda}_\alpha - \overline{\lambda}_\beta)^2 + \left(R_{ij} - \overline{R}_{\gamma\delta}f_i^\gamma f_j^\delta\right)f_k^\alpha f_l^\beta g^{ik}g^{jl}\overline{g}_{\alpha\beta}$$

$$= \sum_{\alpha<\beta} \overline{sec}(\bar{e}_\alpha, \bar{e}_\beta)(\bar{\lambda}_\alpha - \bar{\lambda}_\beta)^2 + g(Ric - f^*\overline{Ric}, f^*\bar{g}) \qquad (3.4)$$

Third, based on the equality (3.4) and the equation (2.2), we formulate the following

**Theorem 1**. *Let $f:(M,g) \to (\overline{M},\bar{g})$ be a harmonic mapping between Riemannian manifolds $(M, g)$ and $(\overline{M},\bar{g})$. Assume that the sectional curvature of the second manifold $(\overline{M},\bar{g})$ is non-negative at every point of $f(M)$ and the first manifold $(M, g)$ is a compact manifold with the Ricci tensor $Ric \geq f^*\overline{Ric}$. Then $f$ is a totally geodesic mapping with constant energy density $e(f)$. Furthermore, if there is at least one point of $M$ at which $Ric > f^*\overline{Ric}$, then $f$ is a constant mapping.*

**Proof**. Let $f:(M,g) \to (\overline{M},\bar{g})$ be a harmonic mapping from a compact Riemannian manifold $(M, g)$ to a Riemannian manifold $(\overline{M},\bar{g})$. We may assume that $M$ is orientable; otherwise, we have only to consider the orientable twofold covering space of $M$. Then by the well known *Green's theorem* (see [14, p. 281]) we obtain from (3.1) the following integral equality

$$\int_M Q(f)\, dVol_g = -\int_M \|Df_*\|^2\, dVol_g. \qquad (3.5)$$

If the inequality $\overline{sec} \geq 0$ is satisfied anywhere on $f(M) \subset \overline{M}$ and the inequality $Ric \geq f^*\overline{Ric}$ is satisfied anywhere on $M$, then $Q(f)$ is non-negative everywhere on $M$. Since our hypothesis implies that the left hand side of (3.5) is non-negative and the right hand side is non-positive, this makes both sides zero. In particular, the left hand side of the above equation gives $\|Df_*\|^2 = 0$, making $f$ a totally geodesic mapping (see [6, p 389]; [7, p. 123]; [31]). On the other hand, the right hand side of the above equation gives $Q(f) = 0$. Then from the equation (3.1) we obtain $\Delta e(f) = 0$ that means $e(f)$ is constant. Next, if the inequalities $\overline{sec} \geq 0$ and $Ric \geq f^*\overline{Ric}$ are satisfied and there is a one point $x$ of $M$ in which $Ric > f^*\overline{Ric}$ then the inequality $\int_M Q(f)\, dVol_g > 0$ is true. This inequality contradicts the equation (3.5). In this case, the harmonic mapping $f$ must be constant. QED.

Consider now the case $\bar{n} = 3$. We know that in dimension three a metric $\bar{g}$ has positive sectional curvature $\overline{sec}$ if and only if $\overline{Ric} < \frac{1}{2}\bar{s}\bar{g}$ for the scalar curvature $\bar{s}$ of $(\overline{M}, \bar{g})$ (see [12]). This proposition follows from the equality (see [34])

$$\overline{sec}(\pi) = \frac{1}{2}\bar{s} - \overline{Ric}(X, X)$$

where $X$ is a unite vector orthogonal to $\pi \subset T_y\overline{M}$ for each point $y \in \overline{M}$. Therefore, if $\bar{n} = 3$ and $\overline{sec} \geq 0$ at each point $f(x) \in \overline{M}$ then the inequality $\overline{Ric} \leq \frac{1}{2}\bar{s}\bar{g}$ is true at each point $f(x) \in \overline{M}$. The above arguments and Theorem 1 allow us to formulate a corollary.

**Corollary 1**. *Let $f : (M, g) \to (\overline{M}, \bar{g})$ be a harmonic mapping between a compact Riemannian manifold (M, g) and a three-dimensional Riemannian manifold $(\overline{M}, \bar{g})$. If the Ricci tensor of $(\overline{M}, \bar{g})$ satisfies the following conditions $f^*\overline{Ric} \leq Ric$ and $\overline{Ric} < \frac{1}{2}\bar{s}\bar{g}$ at each pair of points $x \in M$ and $f(x) \in \overline{M}$, then f is a totally geodesic mapping with constant energy density $e(f)$. Furthermore, if there is at least one point of M at which $f^*\overline{Ric} < Ric$, then f is a constant mapping.*

Next, consider the case $\bar{n} = 1$, which we have excluded from consideration of the above. Recall here that $(M, g)$ is a Riemannian manifold with *quasipositive Ricci curvature* if the Ricci curvature of $(M, g)$ is nonnegative and there is at least one point of $(M, g)$ at which the Ricci curvatures in all directions are positive (see [33]). In this case, we have proved the following proposition (see [28]).

**Corollary 2**. *A compact Riemannian manifold of quasipositive Ricci curvature admits no harmonic submersions onto one-dimensional Riemannian manifolds.*

To conclude this section, we state the result about compact harmonic immersed submanifolds. We assume that $(\overline{M}, \bar{g})$ is a complete non-compact Riemannian manifold with $\overline{sec} \geq 0$ and $\bar{n} \geq 3$. We say that $(\overline{M}, \bar{g})$ has *positive curvature at infinity* if outside a compact set all its sectional curvatures are positive (see [9]). This manifold $(\overline{M}, \bar{g})$ has not compact minimal immersed submanifolds (see also [9]). At the same time it is well known that that an isometric immersion is minimal if and only if it is harmonic (see [7, p. 119]). Thus we have the following:

**Corollary 3**. *Let $f:(M,g)\to(\overline{M},\overline{g})$ be an isometric immersion a compact Riemannian manifold (M, g) into a complete non-compact manifold $(\overline{M},\overline{g})$ with non-negative sectional curvature. If $(\overline{M},\overline{g})$ has positive curvature at infinity, then f is not harmonic.*

## 4. Harmonic contraction mappings between compact Riemannian manifolds

Let $f:(M,g)\to(\overline{M},\overline{g})$ be a harmonic contraction mapping between compact Riemannian manifolds (M, g) and $(\overline{M},\overline{g})$ such that (M, g) has the Ricci tensor $Ric > 0$ at each point $x \in M$ and $(\overline{M},\overline{g})$ has the sectional curvature $\overline{sec} \geq 0$ at each point $f(x) \in \overline{M}$. Recall here that the map $f$ is called contraction if $g \geq f^*\overline{g}$ everywhere on M (see [21]).

For each $x \in M$ let $\lambda(x)$ denotes the smallest eigenvalue for *Ric* and let $\lambda = \inf_{x \in M} \lambda(x)$. On the other hand, for every point $y \in f(M)$ let $\overline{\Lambda}(y)$ denotes the largest eigenvalue for $\overline{Ric}$ and let $\overline{\Lambda} = \sup_{y \in f(M)} \overline{\Lambda}(y)$. In this case, from $\lambda \geq \overline{\Lambda}$ can obtain the inequality $\lambda g \geq \overline{\Lambda} f^*\overline{g}$ because $f$ is a contraction mapping. Hence, we can conclude that $Ric \geq f^*\overline{Ric}$. The last inequality is the second from two "non-existence conditions" of our Theorem 1. Taking the above into account, we can formulate the following corollary from our Theorem 1.

**Corollary 4.** *Let $f:(M,g)\to(\overline{M},\overline{g})$ be a harmonic contraction mapping between compact Riemannian manifolds (M, g) and $(\overline{M},\overline{g})$. Assume that the Ricci tensor of (M, g) is positive definite everywhere on M and the sectional curvature of $(\overline{M},\overline{g})$ is non-negative at every point of $f(M)$. If $\lambda \geq \overline{\Lambda}$, then f is totally geodesic with constant energy density $e(f)$. Furthermore, if $\lambda > \overline{\Lambda}$, then f is a constant mapping.*

In particular, let $\overline{n} = 3$. In this case, if the mapping $f:(M,g)\to(\overline{M},\overline{g})$ is contraction then form the inequality $Ric \geq \frac{1}{2}\,\overline{s}\,g$ we obtain $Ric \geq f^*\overline{Ric}$. Furthermore, if there is at least one point of M at which $Ric > \tfrac{1}{2}\,\bar{s}\,g$ then $Ric > f^*\overline{Ric}$. On the basis of these arguments and our Theorem 1, we can formulate the following corollary.

**Corollary 5.** *Let $f:(M,g)\to(\overline{M},\overline{g})$ be a harmonic contraction mapping between a compact Riemannian manifold (M, g) and a three-dimensional Riemannian manifold $(\overline{M},\overline{g})$. If the Ricci tensor $\overline{Ric}$ and the scalar curvature $\bar{s}$ of $(\overline{M},\overline{g})$ satisfy the following conditions $\overline{Ric}<\tfrac{1}{2}\,\bar{s}\,\overline{g}$ and $Ric\geq\tfrac{1}{2}\,\bar{s}\,g$ at each pair of points $x\in M$ and $f(x)\in\overline{M}$, then f is a totally geodesic mapping with constant energy density $e(f)$. Furthermore, if there is at least one point of M at which $Ric>\tfrac{1}{2}\,\bar{s}\,g$ then f is a constant mapping.*

We consider now a harmonic map $f:(M,g)\to\mathbb{S}^{\bar{n}}(1)$ between a compact Riemannian manifold $(M,g)$ and a unite sphere $\mathbb{S}^{\bar{n}}(1)$ with its standard metric $\overline{g}$. Moreover, we assume that $\bar{n}\geq 2$. In this case, the inequality $\overline{sec}\geq 0$ which is the first "non-existence condition" of our Theorem 1 is satisfied. In turn, for a harmonic contraction mapping $f:(M,g)\to(\overline{M},\overline{g})$ from $Ric\geq(\bar{n}-1)g$ we can obtain the inequality $Ric\geq f^*\overline{Ric}$. Then by our Theorem 1 conclude that $f$ is a totally geodesic mapping. On the other hand, if there is at least one point of M at which $Ric>(\bar{n}-1)g$ we obtain the second "non-existence condition" $Ric>f^*\overline{Ric}$ of the Theorem 1. As a result, we can formulate the following corollary.

**Corollary 5.** *Let $f:(M,g)\to\mathbb{S}^{\bar{n}}(1)$ be a harmonic contraction mapping between a compact Riemannian manifold $(M,g)$ and a unite sphere $\mathbb{S}^{\bar{n}}(1)$ with $\bar{n}\geq 2$. Assume that at every point of $(M,g)$ its Ricci tensor $Ric\geq(\bar{n}-1)g$, then f is a totally geodesic mapping. Furthermore, if there is at least one point of M at which $Ric>(\bar{n}-1)g$ then f is a constant mapping.*

**Remark.** To close this section, we recall that by the well known *Theorem of Bonnet-Myers*, the condition $Ric\geq(n-1)g$ implies $diam(M,g)\leq\pi$ (see [13, p. 194]). Furthermore, if $diam(M,g)=\pi$, then (M, g) is isometric $\mathbb{S}^n(1)$.

## 5. Harmonic mappings of complete Riemannian manifolds to non-negatively curved Riemannian manifolds

Let $f:(M,g)\to(\overline{M},\overline{g})$ be a harmonic mapping between Riemannian manifolds and $(M,g)$ be complete. Assume that the sectional curvature $\overline{sec}\geq 0$ at every point of $f(M)$ then Ricci tensor $\overline{Ric}(X_x,X_x)\geq 0$ at each point $x\in f(M)$ in the direction $X_x\in T_x\overline{M}$ because $\overline{Ric}(X_x,X_x)=\overline{g}(X_x,X_x)\sum_{a=1}^{\overline{n}-1}\overline{sec}(X_x,\overline{e}_a)$ where $\overline{e}_1,...,\overline{e}_{\overline{n}-1}$ an orthonormal basis of $X_x^\perp$. In this case, if we assume that the inequality $Ric\geq f^*\overline{Ric}$ is satisfied anywhere on $M$ then the Ricci tensor $Ric$ is non-negative, and therefore $Q(f)$ is non-negative everywhere on $M$. Schoen and Yau have showed in [22] that $\sqrt{e(f)}$ is subharmonic function on $(M, g)$ if $Q(f)\geq 0$. On the other hand, Yau has proved in his other paper [35] that every non-negative $L^2$-integrable subharmonic function on a complete Riemannian manifold must be constant. Applying this to $\sqrt{e(f)}$, we conclude that $\sqrt{e(f)}$ is a constant if the energy $E(f)=\int_M\|f_*\|^2 dVol_g<+\infty$ (see also [22]). On the other hand, every complete non-compact Riemannian manifold with nonnegative Ricci curvature has infinite volume (see [17]). In our case, we have $Ric\geq 0$ then the volume of $(M, g)$ is infinite. This forces the constant $e(f)$ to be zero and $f$ to be a constant map (see also [22]). We have proved the following result.

**Theorem 2.** *Suppose $f:(M,g)\to(\overline{M},\overline{g})$ is a harmonic mapping with finite energy. Assume that the sectional curvature of $(\overline{M},\overline{g})$ is non-negative at every point of $f(M)$ and $(M, g)$ is a complete non-compact manifold with the Ricci tensor $Ric\geq f^*\overline{Ric}$. Then f is a constant map.*

We know by [5] and [36, p. 247] that harmonic maps $(M, g)$ onto a unite circle $\mathbb{S}^1(1)$ are canonically identified with the harmonic one-forms on $(M,g)$ with integral periods. At the same time, as is well known from the paper [35] that any complete Rie-mannian manifold $(M,g)$ with quasipositive Ricci curvature does not admit harmonic one-forms. Therefore, the following statement is true.

**Corollary 7**. *An arbitrary complete Riemannian manifold with quasipositive Ricci curvature does not admit harmonic maps onto* $\mathbb{S}^1(1)$.

**Remark.** One should compare the result of the Corollary 7 with the result which we stated in the Corollary 5.

## 6. Holomorphic mappings of compact and complete almost semi-Kählerian manifolds

An *almost Hermitian manifold* (*M, g, J*) is a 2*m*-dimensional Riemannian manifold (*M, g*) endowed with almost complex structure *J* which is a smooth section of the bundle $T^*M \otimes TM$ such that $J^2 = -id$ and $g(J,J) = g$ (see [18, p. 147]). In particular, (*M, g, J*) is called *cosymplectic* (see [15, p. 251]) or *almost semi-Kählerian manifold* if $\nabla^* J = 0$ where the operator $\nabla^*$ is formally adjoint to $\nabla$ (see [8]; [10]). The almost semi-Kählerian manifold (*M, g, J*) can be specialized as follows: this manifold is said to be *quasi-Kählerian, nearly-Kählerian* and *Kählerian* if $(\nabla_X J)Y + (\nabla_{JX} JY) = 0$, $(\nabla_X J)X = 0$ and $\nabla J = 0$, respectively (see [10]).

Now, let us consider a smooth mapping $f : (M, g, J) \to (\overline{M}, \overline{g}, \overline{J})$ between almost Hermitian manifolds. We suppose that *f* is not necessarily a diffeomorphism and the manifolds (*M, g, J*) and $(\overline{M}, \overline{g}, \overline{J})$ do not have the same dimension. The mapping *f* is called *holomorphic* if its differential $f_*$ commutes with almost complex structures *J* and $\overline{J}$, i.e. $f_* \circ J = \overline{J} \circ f_*$ (see [15, p. 123]). On the other hand, the mapping *f* is said to be anti-*holomorphic* if $f_* \circ J = -\overline{J} \circ f_*$ (see [18, p. 9]). The basic relation between holomorphic (anti-holomorphic) maps and harmonic maps is given by the following local result due to Lichnerowicz (see [1, p. 252]; [18, p. 9]; [19; p. 175]).

**Theorem 3**. *Any holomorphic map* $f : (M, g, J) \to (\overline{M}, \overline{g}, \overline{J})$ *from an almost semi-Kählerian manifold* (*M, g, J*) *to a quasi-Kählerian manifold* $(\overline{M}, \overline{g}, \overline{J})$ *is harmonic.*

Then using Theorem 2 and Theorem 3, we can state the following corollaries.

**Corollary 8.** *Let* $(\overline{M}, \overline{g}, \overline{J})$ *be a quasi-Kählerian manifold with the sectional curvature* $\overline{sec} \geq 0$ *and* $(M, g, J)$ *be a compact almost semi-Kählerian (quasi-Kählerian, nearly-Kählerian and Kählerian) manifold with the Ricci tensor* $Ric \geq f^*\overline{Ric} \geq 0$. *Then any holomorphic mapping* $f: (M, g, J) \to (\overline{M}, \overline{g}, \overline{J})$ *is totally geodesic. Furthermore, if there is at least one point of M at which* $Ric > f^*\overline{Ric}$, *then f is constant.*

**Corollary 9.** *Let* $(\overline{M}, \overline{g}, \overline{J})$ *be a quasi-Kählerian manifold with the sectional curvature* $\overline{sec} \geq 0$ *and* $(M, g, J)$ *be a complete non-compact almost semi-Kählerian (quasi-Kählerian, nearly-Kählerian and Kählerian) manifold with the Ricci tensor* $Ric \geq f^*\overline{Ric} \geq 0$. *Then any holomorphic mapping* $f: (M, g, J) \to (\overline{M}, \overline{g}, \overline{J})$ *with finite energy is constant.*

## 7. A uniqueness theorem for Ricci tensor

The main point of the papers [4]; [12] and the monograph [2; pp. 140-153] is that in certain circumstances the metric (or at last the connection) is uniquely determined by the Ricci tensor. In reticular, in [4, Corollary 3.3] and [2, Theorem 5.42] the celebrated result was proved: Let $(M, \overline{g})$ be a compact Einstein manifold with the Ricci tensor $\overline{Ric} = \overline{g}$ and the section curvature $\overline{sec} \leq 0$, then an another Riemannian metric $g$ on $M$ with $Ric = \overline{g}$ has the Levi-Civita connection $\nabla$ such that $\nabla = \overline{\nabla}$. At the same time, from Theorem 1 we obtain the corollary which generalizes this result.

**Theorem 3.** *Let* $(M, \overline{g})$ *be a compact Riemannian manifold with the sectional curvature* $\overline{sec} \geq 0$ *and the Ricci tensor* $\overline{Ric} \leq \overline{g}$. *If g is another Riemannian metric on M with* $Ric = \overline{g}$, *then g and* $\overline{g}$ *have the same Levi-Civita connection. Furthermore, if the full holonomy group* $\mathrm{Hol}(\overline{g})$ *is irreducible then* $Ric = \overline{Ric}$.

**Proof.** We assume that the manifold $M$ is equipped with two Riemannian metrics $g$ and $\overline{g}$. If $trace_g T = 0$ for $T = \overline{\nabla} - \nabla$, then the identity mapping from $M$ with metric $g$ to $M$ with metric $\overline{g}$ is harmonic, and $trace_g \overline{g}$ is the harmonic mapping energy density (see [14]). In our situation, we have $Ric = \overline{g} > 0$, then the identity map Id:

$(M, g) \to (M, Ric)$ is harmonic and its mapping energy density $trace_g \bar{g} = s > 0$ for the scalar curvature $s$ of $g$ (see [2, p. 152]; [12, pp. 50-51]). In this case, the Eells-Sampson equation (3.1) has the form (see [4])

$$\tfrac{1}{2} \Delta s = Q(f) + \|T\|^2 \qquad (6.1)$$

where $Q(f) = g^{ik} g^{jl} (\bar{g}_{ij} \bar{g}_{kl} - \bar{R}_{ijkl})$ and $\|T\|^2 = g^{ij} g^{kl} \bar{g}_{pq} T_{ik}^p T_{jl}^q \geq 0$. On the other hand, we have the identity (see [2, p. 436]; [3])

$$(\bar{g}_{ij} \bar{R}_{kl} - \bar{R}_{ijkl}) \varphi^{ik} \varphi^{jl} = \sum_{i<j} \overline{sec}(\bar{e}_i, \bar{e}_j)(\bar{\lambda}_i - \bar{\lambda}_j)^2 \qquad (6.2)$$

where $\varphi$ is any smooth symmetric tensor field such as $\varphi(\bar{e}_i, \bar{e}_j) = \bar{\lambda}_i \delta_{ij}$ for the Kronecker delta $\delta_{ij}$ and an orthonormal basis $\{\bar{e}_1, ..., \bar{e}_n\}$ at any point $x \in M$. Then equation (6.1) can be rewritten in the form

$$\tfrac{1}{2} \Delta s = \sum_{i \neq j} \overline{sec}(\bar{e}_i, \bar{e}_j)(\bar{\lambda}_i - \bar{\lambda}_j)^2 + g^{ik} g^{jl} \bar{g}_{ij}(\bar{g}_{kl} - \bar{R}_{kl}) + \|T\|^2 \qquad (6.3)$$

where $g(\bar{e}_i, \bar{e}_j) = \bar{\lambda}_i \delta_{ij}$. The proof of Corollary 2 is completed by showing that under the stated assumptions the right side of (6.3) is non-negative. Since then $\Delta s \geq 0$ and hence $s$ is a positive subharmonic function on $(M, g)$. If $(M, \bar{g})$ is a compact Riemannian manifold, then using the *Hopf's lemma* (see [15, p. 338]), one can verify that $s = $ const. In this case, from (3.3) we obtain $T = 0$. Therefore, $g$ and $\bar{g}$ have the same Levi-Civita connection. Furthermore, if the full holonomy group Hol($\bar{g}$) of $(M, \bar{g})$ is irreducible, then the metric $g = C \bar{g}$ for some constant $C > 0$ (see [2, pp. 282; 285-287]). In this case, we have the identity $Ric = \overline{Ric}$ because the Ricci tensors $\overline{Ric}$ and $Ric$ of the metrics $\bar{g}$ and $C\bar{g}$, respectively, are equal. QED.

The following two propositions will clarify the situation for vanishing theorems in this theory. In [4] the following vanishing theorem was proved: Let $(M, \bar{g})$ be a compact Riemannian manifold with all sectional curvature less then $(\bar{n} - 1)^{-1}$. Then there is no Riemannian metric $g$ on $M$ such that its Ricci tensor $Ric = \bar{g}$. At the same time, we note that the inequality $\overline{sec} < +1$ implies the following inequality $\overline{Ric} < \bar{g}$. In its turn, in [12] the following non-existence theorem was proved: Let $\bar{g}$ be a metric on a compact manifold $M$ with the sectional curvature $\overline{sec} < +1$, then any metric $g$ does not exist on $M$ such that its Ricci tensor $Ric = \overline{g}$. We also get a non-existence result which complements the above propositions.

**Corollary 10**. *Let $(M, \overline{g})$ be a compact Riemannian manifold with the nonnegative section curvatures and the Ricci tensor $\overline{Ric} < \overline{g}$. Then there does not exist Riemannian metric $g$ on $M$ such that its Ricci tensor $Ric = \overline{g}$.*

**Proof.** Let $(M, \overline{g})$ be a compact Riemannian manifold and the Ricci satisfies the inequality $\overline{Ric} < \overline{g}$. In this case, from the equation (6.3) we obtain the inequality $\Delta s > 0$, which is impossible at a maximum of $s$ we have $\Delta s \leq 0$. QED.

In particular, for $n = 3$ and $\overline{sec} \geq 0$ from (6.3) we obtain

$$\tfrac{1}{2}\,\Delta s \geq \sum_{i<j} \overline{sec}(\overline{e}_i, \overline{e}_j)(\overline{\lambda}_i - \overline{\lambda}_j)^2 + \tfrac{1}{2}\,\|\overline{g}\|^2(2 - \overline{s}) + \|T\|^2 \qquad (6.4)$$

because $\overline{Ric} \leq \tfrac{1}{2}\,\overline{s}\,\overline{g}$. If we suppose that $\overline{s} \leq 2$ then from (6.4) we obtain $\Delta s \geq 0$ and hence $s$ is a positive subharmonic function on $(M, g)$. Furthermore, if $\overline{s} < 2$ then $\Delta s > 0$. Therefore, the Theorem 3 has an impotent corollary.

**Corollary 11**. *Let $(M, \overline{g})$ be a three-dimensional compact Riemannian manifold with the Ricci tensor $\overline{Ric}$ such that $0 \leq \overline{Ric} \leq \tfrac{1}{2}\,\overline{s}\,\overline{g}$ for the scalar curvature $\overline{s} \leq 2$. If $g$ is another Riemannian metric on $M$ with $Ric = \overline{g}$, then $g$ and $\overline{g}$ have the same Levi-Civita connection. Furthermore, if $0 \leq \overline{Ric} \leq \tfrac{1}{2}\,\overline{s}\,\overline{g}$ for $\overline{s} < 2$, then there does not exist Riemannian metric $g$ on $M$ such that its Ricci tensor $Ric = \overline{g}$.*

## References

1. Bair P., Wood J.C.: Harmonic morphisms between Riemannian manifolds, Clarendon Press, Oxford, 2003.
2. Becce, A.L.: Einstein manifolds, Springer-Verlag, Berlin · Heidelberg (1987).
3. Berger, M., Ebin, D.: Some decompositions of the space of symmetric tensors on a Riemannian manifold, Journal of Differential Geometry, **3**(3-4), 379-392 (1969).

4. De Turck, D., Koiso, N.: Uniqueness and non-existence of metrics with prescribed Ricci curvature. Annales de l'Institut Henri Poincare, Section A (N.S), **1**(5), 351-359 (1984).

5. Eells, J., Lemaitre, L.: A report on harmonic maps, Bull. London Mathematical Society, **10**(1), 1-68 (1978).

6. Eells, J., Lemaire, L.: Another report on harmonic maps, Bull. London Mathematical Society, **20**(5), 385-524 (1988).

7. Eells, J., Sampson, J.H.: Harmonic mappings of Riemannian manifolds, American Journal of Mathematics, **86**(1), 109-160 (1964).

8. Friedland, L., Hsiung, Ch.-Ch.: A certain class of almost Hermitian manifolds, Tensor (N.S.), **48**, 252-263 (1989).

9. Galloway G., Rodriguez L.: Intersections of minimal submanifolds, Geometriae Dedicata, **31** (1), 29-42 (1991).

10. Gray, A., Hervella, L.M.: The sixteen classes of almost Hermitian manifolds and their linear invariants, Annali di Matematica Pura ed Applicata, **123**(1), 35-58 (1980).

11. Hamilton R.S.: Three-manifolds with positive Ricci curvature, Journal of Differential Geometry, **17** (2), 255-306 (1982).

12. Hamilton R.S.: The Ricci curvature equation, Lecture notes: Seminar on nonlinear partial differential equations, Mathematical Sciences Research Institute Publications, Berkeley, 47-72 (1983).

13. Jost, J.: Riemannian geometry and geometric analysis, Springer Science and Business Media, Berlin - Heidelberg (2008).

14. Kobayashi, S., Nomizu, K.: Foundations of Differential Geometry, Vol. 1, Interscience Publishers (Wiley), New York (1963).

15. Kobayashi, S., Nomizu, K.: Foundations of Differential Geometry, Vol. 2, Interscience Publishers (Wiley), New York (1969).

16. Lee K.-W., Lee Y.-I.: Mean curvature flow of the graphs of maps between compact manifolds, Transactions of the American Mathematical Society, **363** (11), 5745-5759 (2011).

17. Lichnerowicz, A.: Applications harmoniques et variétés Kähleriennes, Rendiconti dei Seminario Matematico e Fisico di Milano, **39** (1) 186-195 (1969).

18. Pigola, S., Rigoli, M., Setti, A.G.: Vanishing and finiteness results in geometric analysis. A generalization of the Bochner technique, Birkhäuser, Basel (2008).

19. Salamon S.: Harmonic and holomorphic maps, Geometry Seminar "Luigi Bianchi", Pisa 1984, Lecture Notes Mathematics, 164, 161-224 (1985).

20. Sampson, J.H.: Applications of harmonic maps to Kähler geometry, Complex differential geometry and non-linear differential equations, Contemp. Math. **49**, 125-134 (1986).

21. Savas-Halilaj A., Smoczyk K.: Evolution of contractions by mean curvature flow, Mathematische Annalen, **361**(3) 725-740 (2014).

22. Schoen, R., Yau, S.T.: Harmonic maps and topology of stable hypersurfaces and manifolds with non-negative Ricci curvature, Commenttarii Mathematici Helvetici, **51**(1), 333-341 (1976).

23. Sealey, H.C.J.: Some properties of harmonic mappings, PhD thesis, University of Warwick (1980).

24. Siu, Y.T.: The complex-analyticity of harmonic maps and the strong rigidity of compact Kähler manifolds, Annals of Mathematics, **12**, 73-111 (1980).

25. Stepanov, S.E.: Vanishing theorems in affine, Riemannian, and Lorentz geometries, Journal of Mathematical Sciences (NY), **141**(1), 929-964 (2007).

26. Stepanov, S.E.: On the global theory of some classes of mappings, Annals of Global Analysis and Geometry, **13**(3), 239-249 (1995).

27. Stepanov S.E.: Geometry of infinitesimal harmonic transformations, Annals of Global Analysis and Geometry, **24**(3), 291-299 (2003).

28. Stepanov S.E.: O($n$) × O($m - n$)-structures on $m$-dimensional manifolds and submersions of Riemannian manifolds, St. Petersburg Mathematical Journal, **7** (6), 1005-1016 (1996).

29. Stepanov S., Tsyganok I.: Vanishing theorems for projective and harmonic mappings, Journal of Geometry, **106** (3), 640-641 (2015).

30. Tsui M.-P., Wang M.-T.: Mean curvature flows and isotopy of maps between spheres, Communications on Pure and Applied Mathematics, **57** (8), 1110-1126 (2004).

31. Vilms, J.: Totally geodesic maps, Journal Differential Geometry, **4**(1), 73-79 (1970).

32. Wu, H.H.: The Bochner technique in differential geometry, Mathematical Reports, Harwood Acad. Publ. (1987).

33. Wu H.: A remark on the Bochner technique in differential geometry, Proceedings of the American Mathematical Society, **78** (2), 403-408 (1980).

34. Wolfson J., Schmidt B.: Three-manifolds with constant vector curvature, Indiana University Mathematics Journal, **63** (6), 1757-1783 (2014).

35. Yau, S.T.: Some function-theoretic properties of complete Riemannian manifold and their applications to geometry, Indiana University Mathematical Journal, **25**(7), 659-679 (1976).

36. Yau, S.T., Jost J.: Harmonic maps and superrigidity, Proceedings of Symposia in Pure Mathematics, **54** (1), 245-326 (1993).